\documentclass[12pt]{iopart}
\usepackage{graphicx}
\usepackage{xcolor}
\usepackage{epstopdf}

\usepackage{bm}

\usepackage{amsfonts}
\usepackage{amssymb}
\usepackage{amsthm}
\usepackage{bbm}

\usepackage{verbatim}
\usepackage[textsize=footnotesize]{todonotes}

\newtheorem{Theorem}{Theorem}
\newtheorem*{theorem}{Theorem}
\newtheorem{Lemma}[Theorem]{Lemma}
\newtheorem{Propo}[Theorem]{Proposition}
\newtheorem{Coro}[Theorem]{Corollary}
\newtheorem{Definition}[Theorem]{Definition}
\newtheorem{Remark}[Theorem]{Remark}

\newcommand{\abs}[1]{{\left|#1\right|}}
\newcommand{\norm}[1]{{\left\|#1\right\|}}
\renewcommand{\Re}{\mathrm{Re}}
\renewcommand{\d}{\;{\rm d}} 
\renewcommand{\epsilon}{\varepsilon}
\newcommand{\fa}          {\quad \mbox{for all } \,}
\newcommand{\set}[1]      {\{#1\}}

\newcommand{\R}{\mathbb{R}}
\newcommand{\1}{\mathbbm{1}}

\usepackage{color}

\begin{document}

\title[Coupling functions allowing persistent synchronisation]
{Towards a general theory for coupling functions allowing persistent synchronisation}

\author{Tiago Pereira$^1$, Jaap Eldering$^1$, Martin Rasmussen$^{1}$, and Alexei Veneziani$^2$}
\address{$^1$Department of Mathematics, Imperial College London, London SW7 2AZ, UK \\
$^2$Centro de Matem\'{a}tica, Computa\c{c}\~{a}o e Cogni\c{c}\~{a}o, UFABC, Santo Andr\'{e}, Brazil }
\ead{tiago.pereira@imperial.ac.uk, j.eldering@imperial.ac.uk, m.rasmussen@imperial.ac.uk}

\begin{abstract}

We study synchronisation properties of networks of coupled dynamical systems with interaction akin to diffusion. We assume that the isolated node dynamics possesses a forward invariant set on which it has a bounded Jacobian, then we characterise a class of coupling functions that allows for uniformly stable synchronisation in connected complex networks --- in the sense that there is an open neighbourhood of the initial conditions that is uniformly attracted towards synchronisation. Moreover, this stable synchronisation persists under perturbations to non-identical node dynamics. We illustrate the theory with numerical examples and conclude with a discussion on embedding these results in a more general framework of spectral dichotomies.
\end{abstract}

\maketitle


\section{Introduction}

Network synchronisation is observed to occur in a broad range of applications in physics~\cite{Wiesenfeld_98_1}, neuroscience~\cite{Bullmore_09_1,Gregoriou_09_1,Singer_99_1,Milton_03_1}, and ecology~\cite{Earn_00_1}. During the last fifty years, empirical studies of real complex systems have led to a deep understanding of the structure of networks~\cite{Newman_10_1,Albert_02_1}, and the interaction properties between oscillators, that is, the coupling function~\cite{Kuramoto_84_1,Stankovski_12_1,Wu_07_1}.

The stability of network synchronisation is a balance between the isolated dynamics and the coupling function. Past research suggests that in networks of identical oscillators with interaction akin to diffusion, under mild conditions on the isolated dynamics, the coupling function dictates the synchronisation properties of the network~\cite{Pecora_98_1,Li_04_2,Pereira_10_1,Orosz_09_1,Wu_07_1}. However, it still remains an open problem to describe the class of coupling functions that lead the network to persistent synchronisation.

Our work contributes to the development a general theory for coupling functions that allow for persistent synchronisation for a connected complex network. The coupling functions under consideration appear in a variety of synchronisation models on networks (such as the Kuramoto models~\cite{Kuramoto_84_1} and its extensions~\cite{Acebron_05_1,Blasius_05_1,Petkoski_12_1}).

More precisely, we consider the dynamics of a network of $n$ identical elements with interaction akin to diffusion, described by
\begin{eqnarray}
  \dot{x}_i = f(t, x_i) + \alpha \sum_{j=1}^n W_{ij} h(x_j  - x_i)\,,
  \label{eq:md1}
\end{eqnarray}
where $\alpha$ is the overall coupling strength, and the matrix $W = (W_{ij})_{i,j\in\{1,\dots,n\}}$ describes the interaction structure of the network, i.e.~$W_{ij}$ measures the strength of interaction between the nodes $i$ and $j$. The function $f: \mathbb{R} \times \mathbb{R}^m \rightarrow \mathbb{R}^m$ describes the isolated node dynamics, and the coupling function $h : \mathbb{R}^m \rightarrow \mathbb{R}^m$ describes the diffusion-like interaction between nodes. We make the following two assumptions for these functions.

\medskip

\noindent \textbf{Assumption A1.} \textsl{The function $f$ is continuous, and there exists an inflowing invariant open ball $U \subset \mathbb{R}^m$ such that $f$ is continuously differentiable in $U$ with
\begin{displaymath}
  \| D_2 f(t,x)\| \le \varrho \fa t \in \mathbb{R} \mbox{ and } x \in U
\end{displaymath}
for some $\varrho>0$.}

\medskip

For instance, the Lorenz system has a bounded inflowing invariant ball, see Subsection~\ref{Lorenz}. In general, smooth nonlinear systems with compact attractors satisfy Assumption~A1. This assumption will be generalised in Section~\ref{auxiliary} to include also noncompact sets $U$.

\medskip

\noindent \textbf{Assumption A2.} \textsl{The coupling function $h$ is continuously differentiable with $h(0) = 0$. We define $\Gamma := Dh(0)$ and denote the (complex) eigenvalues of $\Gamma$ by $\beta_i$, $i\in\{1,\dots,m\}$.}

\medskip

The network structure plays a central role for the synchronisation properties. We consider the intensity of the $i$-th node $V_i = \sum_{j=1}^n W_{ij}$, and define the positive definite matrix $V := \mathrm{diag}(V_1,\dots,V_n)$.  Then the so-called \emph{Laplacian} reads as
\[
  L = V - W\,.
\]
Let $\lambda_i$, $i\in\{1,\dots,n\}$, denote the eigenvalues of $L$. Note that $\lambda_1=0$ is an eigenvalue with eigenvector $\frac{1}{\sqrt{n}}(1, \dots , 1)$. The multiplicity of this eigenvalue equals the number of connected components of the network.

The following assumption incorporates the coupling and structural network properties.

\medskip

\noindent \textbf{Assumption A3.} \textsl{We suppose that
\[
\gamma := \min_{{ 2\le i \le n} \atop  {1\le j \le m}} \Re(\lambda_i \beta_j ) > 0\,,
\]
where $\Re(z)$ denotes the real part of a complex number $z$.}

\medskip

The dynamics of such a diffusive model can be intricate. Indeed, even if the isolated dynamics possesses a globally stable fixed point, the diffusive coupling can lead to instability of the fixed point and the system can exhibit an oscillatory behaviour~\cite{Pogromsky_99_1}.

Note that due to the diffusive nature of the coupling, if all oscillators start with the same initial condition,
then the coupling term vanishes identically. This ensures that the globally synchronised state $x_1(t)=x_2(t)= \dots = x_n(t) = s(t)$ is an invariant state for all coupling strengths $\alpha$ and all choices of coupling functions $h$. That is, the \emph{diagonal manifold}
\begin{equation*}
  M := \big\{ x_i \in \R^{m} \mbox{ for } i \in \{ 1, \cdots , n \} : x_1=\cdots=x_n \big\}
\end{equation*}
is invariant, and we call the subset
\begin{equation}\label{eq:sync-manifold}
  S := \big\{ x_i \in U \subset \R^{m} \mbox{ for } i \in \{ 1, \cdots , n \} : x_1=\cdots=x_n \big\} \subset M
\end{equation}
the \emph{synchronisation manifold}. The main result of this paper is a proof that under the general conditions given above and $\alpha$ sufficiently large, the synchronisation manifold $S$ is uniformly exponentially stable.

\begin{Theorem}[synchronisation]\label{MainThm}
Consider the network of diffusively coupled equations~(\ref{eq:md1}) satisfying~A1--A3. Then there exists a $\rho = \rho(f,\Gamma)$ such that for all coupling strengths
\[
  \alpha > \frac{\rho}{\gamma} \,,
\]
the network is \emph{locally uniformly synchronised}. This means that there exist a $\delta>0$ and a $C=C(L,\Gamma)>0$ such that if $x_i(t_0) \in U$ and $\| x_i (t_0) - x_j (t_0)\| \le \delta$ for any $i,j\in\{1,\dots,n\}$, then
\begin{equation}\label{eq:thm-exp-rate}
  \| x_i (t) - x_j (t)\| \le C e^{-(\alpha \gamma - \rho)(t-t_0)} \| x_i (t_0) - x_j (t_0)\| \quad \mbox{for all } t\ge t_0\,.
\end{equation}
\end{Theorem}

Hence, the synchronisation manifold is locally uniformly exponentially attractive. The constant $\rho$ depends on
the bounds on the Jacobian $D_2 f$ as set out in Assumption~A1 and on the conditional number of the matrix $\Gamma$ (see~(\ref{rhovarrho}) in case $\Gamma$ is diagonalisable). In the case that the Laplacian $L$ and $\Gamma$ are diagonalisable, $C$ depends on the conditional number of the similarity transformation that diagonalises these matrices (see Lemma~\ref{lem:NetDiag} for details), so loosely speaking, it depends on how well the eigenvectors of $L$ and $\Gamma$ are orthogonal. If $L$ and $\Gamma$ are non-diagonalisable, then $C$ is related to conditional numbers as well, see the proof of Lemma~\ref{lem:Nondiag}
for details. The size of $\delta$ can be estimated explicitly if more concrete details about the system are known, see also  Remark~\ref{rem:delta-est} on page~\pageref{rem:delta-est}.

Our second main result shows that synchronisation is persistent under perturbation of the isolated nodes. Thereto, consider a network of non-identical nodes described by
\begin{eqnarray}\label{eq:perturbed}
\dot{x}_i = f_i(t, x_i) + {\alpha} \sum_{j=1}^n W_{ij} h(x_j  - x_i),
\end{eqnarray}
where $f_i(t, x_i) = f(t, x_i) +g_i(t, x_i)$. Note that in this case, the synchronisation manifold $S$ is no longer invariant. We show in this paper that for small perturbations functions $g_i$, $i\in\{1,\dots,n\}$, the synchronisation manifold is stable in the sense that orbits starting near the synchronisation manifold $S$ remain in a neighbourhood of $S$.

\begin{Theorem}[persistence]\label{PertThm}
  Consider the perturbed network~(\ref{eq:perturbed}) of diffusively coupled equations fulfilling Assumptions A1--A3, and suppose that
  \[
    \alpha > \frac{\rho}{\gamma}
  \]
  as in Theorem~\ref{MainThm}. Then there exist $\delta>0$, $C>0$ and $\epsilon_g>0$ such that for all $\epsilon_0$-perturbations satisfying
  \begin{equation}\label{eq:pert-size}
    \norm{g_i(t, x)} \le \epsilon_0 \le \epsilon_g \quad \mbox{for all } t\in\R\,,\, x\in U\mbox{ and } i\in\{1,\dots, n\}
  \end{equation}
  and initial conditions satisfying $\norm{x_i(t_0) - x_j(t_0)} \le \delta$ for any $i,j\in\{1,\dots,n\}$, the estimate
  \begin{equation}\label{eq:pert-bound}
    \hspace{-1cm}\norm{x_i(t) - x_j(t)}
    \le C e^{-(\alpha \gamma - \rho )(t-t_0)} \norm{x_i(t_0) - x_j(t_0)}
       +\frac{C\epsilon_0}{\alpha \gamma  - \rho}
    \quad \mbox{for all } t \ge t_0
  \end{equation}
  holds.
\end{Theorem}

Note that the additional term $C\epsilon_0 / (\alpha \gamma  - \rho)$ can be made small either by controlling the perturbation size $\epsilon_0$ or by increasing $\alpha\gamma$. This provides control of the network coherence in terms of the network properties and coupling strength.

If the Laplacian $L$ is symmetric (i.e.~the systems are mutually coupled), its spectrum is real and can be ordered as $0 = \lambda_1 < \lambda_2 \le \lambda_3 \le \dots \le \lambda_n$. Moreover, consider $\beta := \min_{i \in \{1,\dots, m\}} \Re \beta_i$, and note that this implies
\[
  \gamma = \beta \lambda_2\,.
\]

The following corollary to the above persistence result then shows that the enhancement of coherence in the network in terms of network connectivity depends on the spectral gap $\lambda_2$.

\begin{Coro}[synchronisation error]\label{corofpersistence}
Consider the perturbed network~(\ref{eq:perturbed}) with symmetric Laplacian $L$ and the average synchronisation error
\[
e_s(t) = \displaystyle{ \frac{1}{n(n-1)} \sum_{i,j=1}^n \| x_i(t) - x_j(t) \|} \fa t\ge t_0\,,
\]
where the initial conditions $x_i(t_0)$, $i\in\set{1,\dots,n}$, are chosen as in Theorem~\ref{PertThm}. Then whenever $\alpha\gamma = \alpha\beta\lambda_2 > \rho$, one has
\[
  \limsup_{t\to\infty} e_s(t)  \le K \frac{\epsilon_0}{\alpha \beta \lambda_2 - \rho} \,,
\]
where $K=K(\Gamma)$ is independent of the network size.
\end{Coro}

This corollary has excellent agreement with recent numerical simulations
for the synchronisation transition in complex networks of mutually
coupled non-identical oscillators~\cite{Pereira_13_1}.

The paper is organised as follows. In Section~\ref{discussion}, we discuss
our assumptions, ideas of the proofs as well as
how our results relates to previous contributions. In Section~\ref{illustrations},
we illustrate our main synchronisation result with a nonautonomous linear system and a coupled Lorenz system. Section~\ref{notation} provides fundamental results on nonautonomous linear differential equations. In Section~\ref{auxiliary}, we provide auxiliary results to prove our main theorems in Sections~\ref{sec:sync} and~\ref{sec_pers}. Finally, in Section~\ref{secgeneral}, we discuss how to generalise this theory using the dichotomy spectrum and normal hyperbolicity.

\medskip

\noindent \textbf{Notation.} We endow the vector space $\mathbb{R}^m$ with the Euclidean norm $\| x \| = \sqrt{\sum_{i=1}^m |x_i  |^2}$ and the associated Euclidean inner product. In addition, we equip the vector space $(\R^m)^n=\R^{nm}$ with the norm
\begin{equation}\label{eq:max-eucl-norm}
  \norm{(x_1,\dots,x_n)} := \max_{i=1,\dots, n} \norm{x_i} \qquad\textrm{where } x_i \in \R^m\,.
\end{equation}
Note that linear operators on the above spaces will be equipped with the induced operator norm. For a given invertible matrix $A\in\mathbb{R}^{d\times d}$, the \emph{conditional number} is defined by  $\kappa(A) = \| A \| \| A^{-1} \|$. Note that the conditional number depends on the underlying operator norm. Finally, the symbol $I_d$ stands for the identity matrix in $\R^d$.

\section{Discussion of the main results}\label{discussion}

This section is devoted to relating our results to the state of the art and to explaining the assumptions and the central ideas of the proofs.

\subsection{State of the art}

Recent research on synchronisation has focused on the role of the coupling function for the stability of network synchronisation. Notably, Pecora and collaborators have developed so-called \emph{master stability functions} to estimate Lyapunov exponents corresponding to the transversal directions of the synchronisation manifold~\cite{Pecora_98_1,Huang_09_1}. In contrast to this approach, we estimate the contraction rate by dichotomy techniques. Our results show that the synchronisation state is locally stable and persistent, and thus stable under small perturbations. This means that the phenomenon of bubbling~\cite{Ashwin_94_1} and riddling~\cite{Heagy_94_1} (which leads to synchronisation loss) will not be observed under our conditions, in contrast to the master stability function approach.

Another aspect of our results is that the synchronisation properties do not depend on diagonalisation properties of the Laplacian. Recently, the master stability function has been extended to include non-diagonalisable Laplacians~\cite{Nishikawa_06_1}. However, these results do not guarantee that an open neighbourhood of the synchronisation manifold will be attracted by the synchronisation manifold, nor do they imply persistence of the synchronisation. In our set-up, these properties follow naturally by means of roughness of exponential dichotomies, which is relevant in applications that are subjected to noise and external influences. Note that the master stability function approach is applicable to a broader class of coupling functions than the ones we consider, but our approach is constructive and making use of further dichotomy techniques and normal hyperbolicity our results can be generalised further, as discussed later in Section~\ref{secgeneral}.

In addition, Pogromsky and Nijmeijer~\cite{Pogromsky_01_1} use control techniques to show that if the coupling function is linear and given by a symmetric positive definite matrix, then the synchronisation manifold is globally asymptotically stable for connected networks. Likewise, Belykh, Belykh and Hasler~\cite{Belykh_04_1} develop a connection graph stability method to obtain global synchronisation for the network, by assuming the existence of a quadratic Lyapunov function associated with the isolated system. In this article, we tackle only local stability properties, but we consider a more general class of coupling functions. However, under additional conditions on the dynamics and coupling functions, it is possible to prove global stability with the techniques we have developed by applying the mean value theorem instead of using Taylor expansions of the vector field.

\subsection{The assumptions}

Our main assumptions are natural and fulfilled by a large class of systems. Assumption A1 concerns the existence of solutions and the boundedness of the Jacobian. Assumption A2 makes it possible to characterise the stability of synchronisation by the linearisation of $h$. Assumption A3 guarantees that the eigenvalues of the tensor $L \otimes \Gamma$ have real part bounded away from zero (except for the trivial eigenvalue).

These hypotheses basically imply that with a finite value of $\alpha$, we are able to damp all the instabilities of the vector field and obtain a stable synchronisation state. If for example, Assumption~A3 is dropped, $\gamma$ may become negative and synchronisation may no longer be possible.

We illustrate the relevance of Assumption~A3 with the following example.
Consider the isolated dynamics $f: \mathbb{R}^2 \rightarrow \mathbb{R}^2$ given by $f(x) = - \epsilon x$.
Moreover, consider three coupled systems
\begin{eqnarray}
\dot{x}_1 &=& f(x_1) + 2 \alpha \Gamma(x_2 - x_1) +\alpha  \Gamma (x_3 - x_1), \nonumber \\
\dot{x}_2 &=& f(x_2) + 2 \alpha \Gamma(x_3 - x_2) \nonumber \\
\dot{x}_3 &=& f(x_3) +   \alpha \Gamma(x_1 - x_3) \nonumber
\end{eqnarray}
with
\[
\Gamma =
\left(
\begin{array}{cc}
2 & 1 \\
-17 & 0
\end{array}
\right)
\quad\mbox{and note that}\quad
L =
\left(
\begin{array}{ccc}
3 & -2 & -1 \\
0 & 2 & -2 \\
-1 & 0 & 1
\end{array}
\right)\,.
\]
The eigenvalues of $L$ are $\lambda_1 = 0$, $\lambda_2 = 3 + i$ and $\lambda_3 = 3 - i$ and
the eigenvalues of $\Gamma$ are $\beta_1 = 1 + 4i $ and $\beta_2  = 1 - 4 i$.
Hence,
\[
\gamma = -1,
\]
and although the isolated dynamics has a stable trivial fixed point, for any $\alpha > \epsilon$
the origin is unstable and there are trajectories of the coupled systems that escape any compact set.
This shows that breaking condition A3 can have severe effects on the dynamics of the coupled systems.

Assumption~A3 has not been considered in the literature to our best knowledge. In the following, we rephrase this condition in the following two special cases:
\begin{itemize}
  \item[(i)]
  \emph{The spectrum of\/ $\Gamma$ is positive.} If $\Gamma$ has a spectrum consisting of only real, positive eigenvalues, then A3 has a representation in terms of the Laplacian. In this case, this condition reads as
  \[
    \Re(\lambda_i) > 0 \quad \mbox{for all } i \not= 1\,,
  \]
  since the Laplacian always has a zero eigenvalue. If the network is connected, this eigenvalue is simple, and by virtue of the disk theorem, a sufficient condition for all other eigenvalues to have positive real part is \emph{positive interaction strength}, i.e.~$W_{ij}>0$ whenever $i$ is connected to $j$, and zero otherwise.
  \item[(ii)]
  \emph{The Laplacian is symmetric.} This is the most studied case in the literature. Assume that the network is connected. Since the spectrum of the Laplacian is real, Assumption~A3 requires that the real part of the spectrum of $\Gamma$ is positive and that the spectrum of the Laplacian is positive apart from the single zero eigenvalue
(or alternatively, that the spectra of $\Gamma$ and the Laplacian are both negative, but note that this is non-physical).
\end{itemize}

\subsection{Ideas of the proofs}

The proofs of our main results rely on identifying the synchronisation problem with a corresponding fixed point problem. We first concentrate on the case of diagonalisable Laplacians, where diagonal dominance (Proposition~\ref{prop:uniCont}) can be used to show that the synchronised state is uniformly asymptotically stable. To obtain the claim for general coupling functions, we make use of the roughness property associated with the equilibrium point (Theorem~\ref{roughness}). The main aspect here is to approximate the coupling function by a diagonalisable one while keeping control of the contraction rates. Finally, the proof for general Laplacians follows from the fact that the set of diagonalisable Laplacians is dense in the space of Laplacians. From these results and the roughness property the main claim follows.

\section{Illustrations}\label{illustrations}

Before proving the two main results of this paper, two examples are discussed.

\subsection{Nonautonomous Linear Equations}\label{LinEq}

Consider the nonautonomous linear equation
\begin{equation}\label{A}
  \dot x = A(t) x
\end{equation}
where
\[
\hspace{-2.5cm}A(t) =
\left(
\begin{array}{cc}
-1 -9 \cos^2(6t) + 12 \sin (6t) \cos(6t) & 12 \cos^2(6t) + 9 \sin (6t) \cos(6t) \\
-12 \sin^2 (6t) + 9 \sin(t) \cos(6t)  & -1 -9 \sin^2(6t) - 12 \sin(6t) \cos (6t)
\end{array}
\right)\,.
\]
This is a prototypical example where the eigenvalues of the time-dependent matrices do not characterise the stability of a nonautonomous linear system. Indeed, the eigenvalues of $A(t)$ are $-1$ and $-10$, independent of $t\in\R$, and a direct computation shows that
\[
x(t) =
\left(
\begin{array}{c}
e^{2t}( \cos(6t) + 2\sin(6t)) + 2 e^{-13 t} (2 \cos(6t) - \sin(6t)) \\
e^{2t}( \cos(6t) - 2\sin(6t)) + 2 e^{-13 t} (2 \cos(6t) - \sin(6t)) \\
\end{array}
\right)
\]
is a solution of the system, which does not converge to $0$ as $t\to\infty$.

Consider now two diffusively coupled systems
\begin{eqnarray*}
  \dot{x}_1 &=& A(t) x_1 + \alpha \Gamma (x_2 - x_1)\,, \\
  \dot{x}_2 &=& A(t) x_2 + \alpha \Gamma (x_1 - x_2)\,,
\end{eqnarray*}
where $\Gamma$ is a real $2\times 2$ matrix. Theorem~\ref{MainThm} yields that it is possible to synchronise these two systems for any coupling matrix with $\beta(\Gamma)>0$. Consider the coupling matrix
\[
\Gamma =
\left(
\begin{array}{cc}
\beta &  1 \\
0 & \beta
\end{array}
\right)\,.
\]
$\Gamma$ is in its Jordan form and non-diagonalisable. The transformation $y = x_1 - x_2$ leads to
\begin{eqnarray}
  \dot{y} &=& \big(A(t) - 2 \alpha \Gamma \big)y\,.
\label{eq:linXi}
\end{eqnarray}
Our main result shows that the trivial solution of~(\ref{eq:linXi}) is stable if $\alpha$ is large enough.

\begin{figure}[htbp] 
   \centering
   \includegraphics[width=3in]{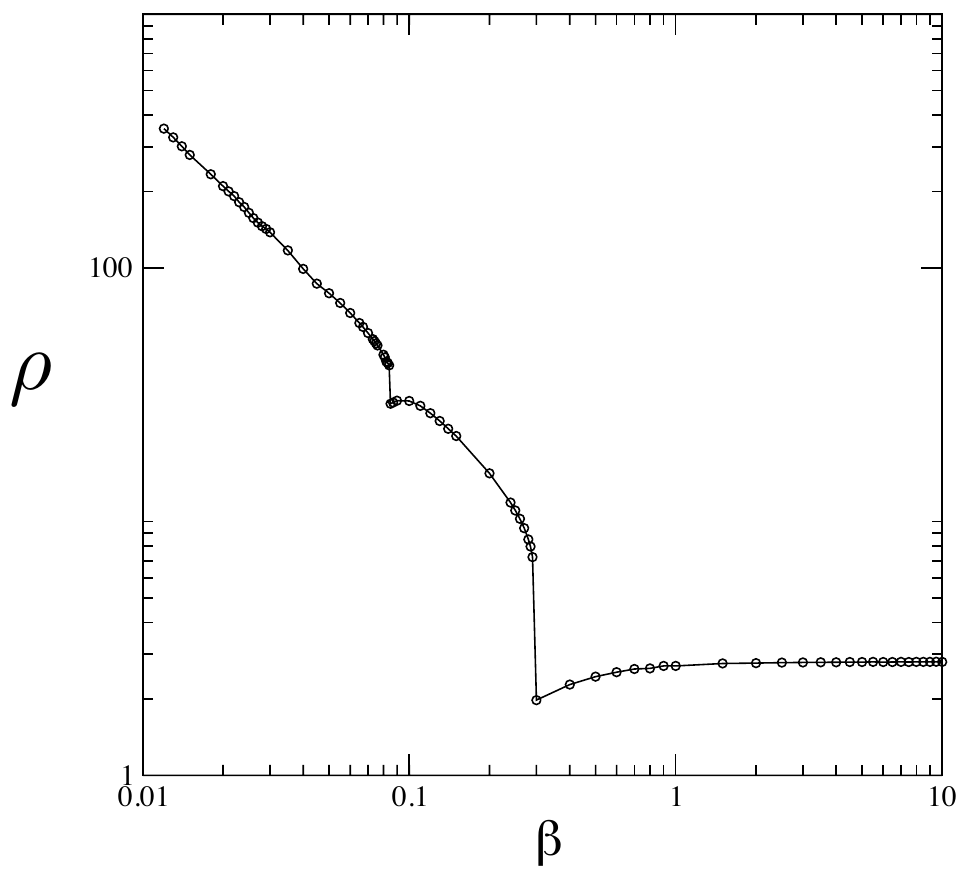}
   \caption{$\rho=\rho(f,\Gamma)$ as a function of $\beta$ in a log--log scale for a fixed $f$ given by Eq.~(\ref{A}). For small $\beta$ the slope is $-1$ in good approximation.}
   \label{fig:LambdaAlpha}
\end{figure}

We have integrated~(\ref{eq:linXi}) using a sixth order Runge--Kutta method with step size $0.001$. We have computed the critical coupling value $\alpha_c$ as a function of $\beta$, such that the trivial solution of Eq.~(\ref{eq:linXi}) is stable.
In Figure~(\ref{fig:LambdaAlpha}) we plotted the corresponding critical value $\rho_c = \beta \alpha_c$. Hence, we are able to analyse the dependence of $\rho$ on $f$ and $\Gamma$.
The behaviour of $\rho$ appears to be intricate. For large $\beta$, we obtain that $\rho$ tends to a constant, however, as we decrease $\beta$, various changes in the behaviour can be observed.
Although the problem is linear, the critical coupling strength depends nonlinearly on the parameter $\beta$. We analyse this dependence
in more details in Section~\ref{Aprho}

\subsection{The Lorenz system}\label{Lorenz}

Using the notation $x =(u, v, w)$, the Lorenz vector field is given by
\[
f(x) =
\left(
\begin{array}{c}
\sigma ( v - u )  \\ u (r - w) - v \\ -b w + u v
\end{array}
\right)\,,
\]
where we choose the classical parameter values $\sigma=10$, $r = 28$ and  $b = \frac{8}{3}$.
All trajectories of the Lorenz system enter a compact set eventually and exist globally forward in time for this reason. Moreover,
they accumulate in a neighbourhood of a chaotic attractor~\cite{Viana_00_1}.

\begin{figure}[htbp] 
   \centering
   \includegraphics[width=4in]{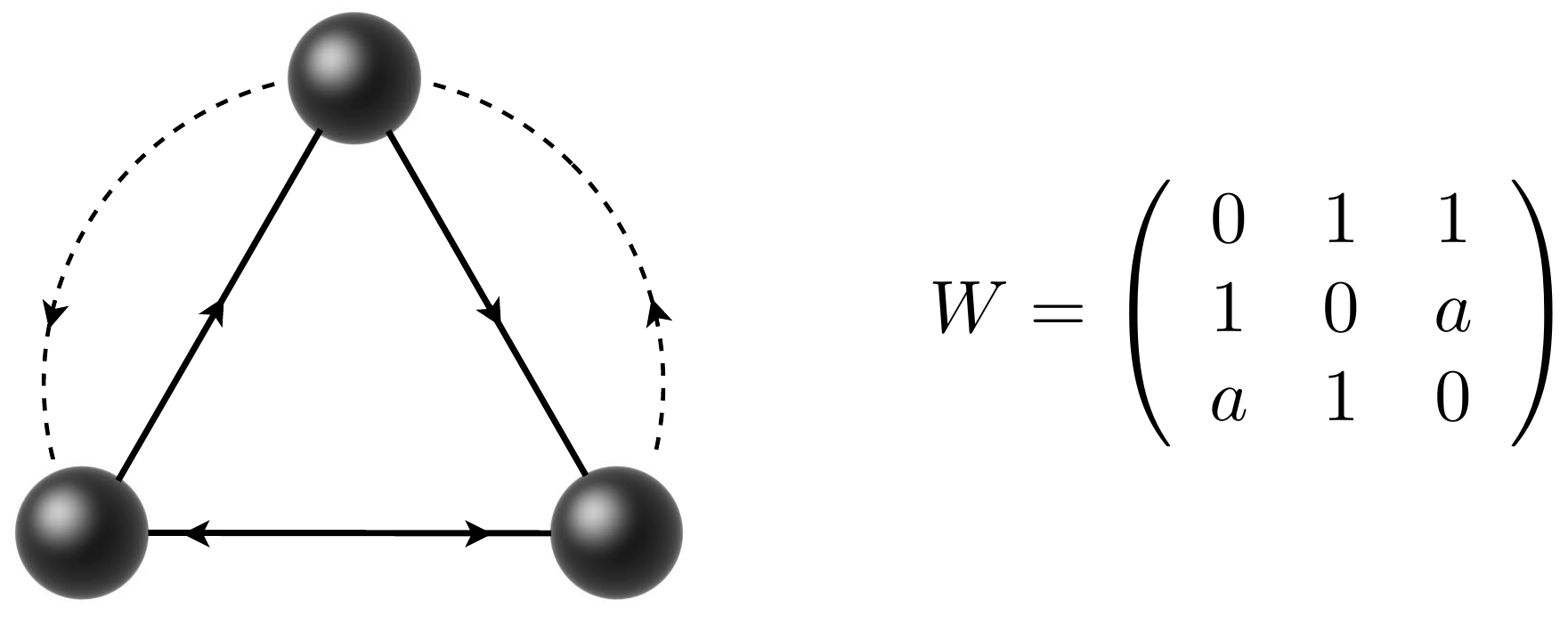}
   \caption{The network and its weight matrix. The matrix $L = V - W$ is non-diagonalisable  for every $a\not=1$; here we choose $a=\frac{1}{3}$.}
   \label{fig:NetL}
\end{figure}

Consider the network of three coupled Lorenz systems
\begin{equation}
  \dot {x}_i = f(x_i) + \alpha \sum_{j=1}^3 W_{ij} H (x_j - x_i)\,,
\label{eq:3Lor}
\end{equation}
where the interaction matrix $W$ is given as in Figure~\ref{fig:NetL}.

\begin{figure}[htbp] 
   \centering
   \includegraphics[width=6in]{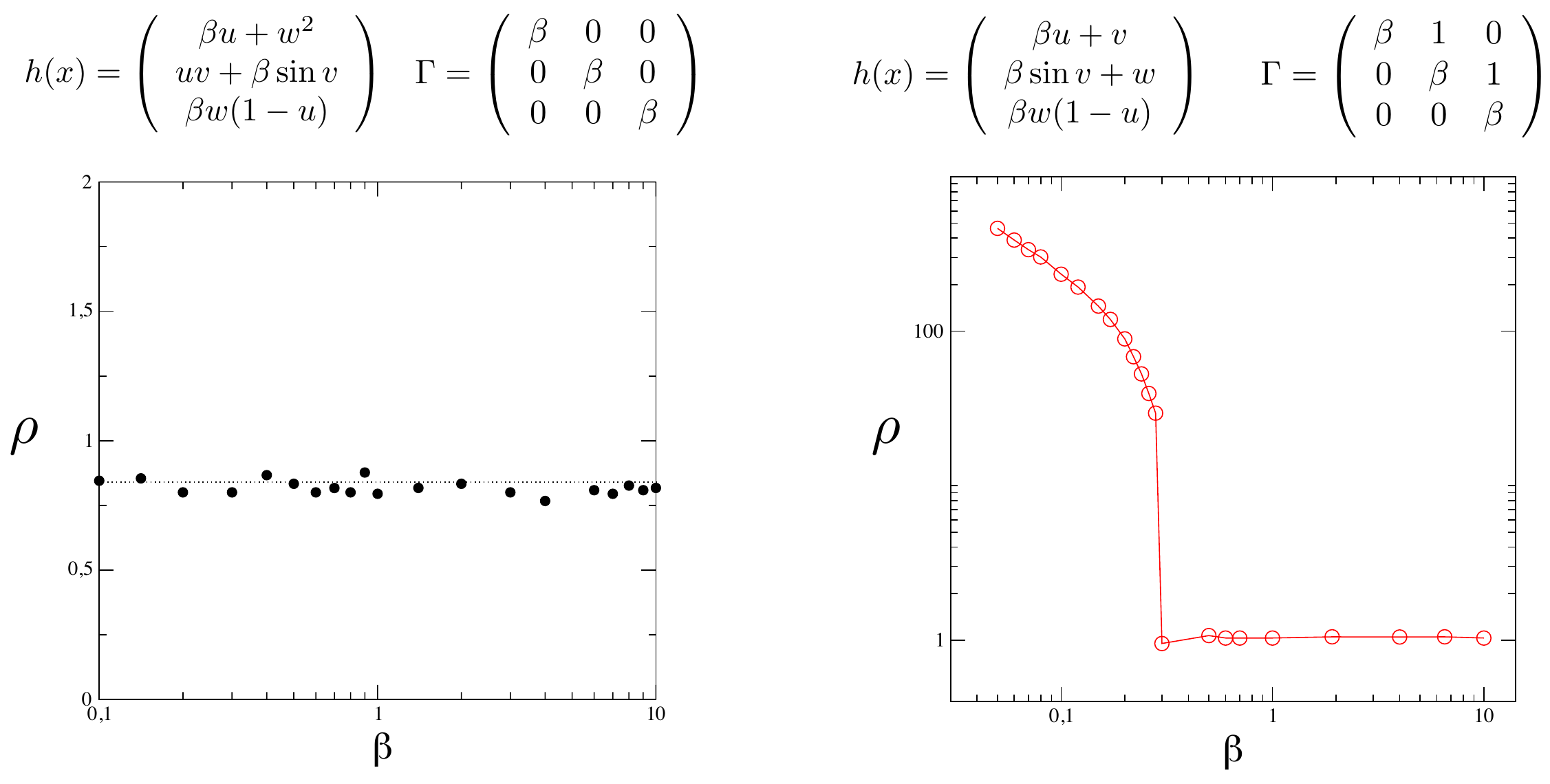}
   \caption{Simulation results for $\rho$ for the
   two coupling functions. For the first case, see left side, $\Gamma = \beta I$
   is positive definite for $\beta>0$, and the behaviour of $\rho$ does not depend significantly
   on $\beta$. For the second case,  $\Gamma$ is a Jordan block with eigenvalues equal to $\beta$.
   In this situation, for large values of $\beta$, the critical coupling $\rho$ appears independent
   of $\beta$, as opposed to the small values of $\beta$. In that case, the critical coupling scales
   as $\rho \propto \beta^{-1}$.}
   \label{fig:Results_Lorenz}
\end{figure}

We use two different nonlinear coupling functions; for the first, the associated matrix $\Gamma$ is positive definite, whereas for the second, $\Gamma$ is a Jordan block. The specific forms of the coupling functions can be seen in Figure~\ref{fig:Results_Lorenz}. We have integrated~(\ref{eq:3Lor}) using a sixth order Runge--Kutta method with step size $0.0001$ and computed the critical coupling $\alpha_c$ as a function of $\beta$, and then plotted the value $\rho_c = \alpha_c \beta$ (see
Figure~\ref{fig:Results_Lorenz}). The behaviour of $\rho$ depends in an essential way on $\Gamma$. This behaviour is further discussed in Section~\ref{Aprho}.

\section{Nonautonomous linear differential equations}\label{notation}

Consider the $m$-dimensional linear differential equation
\begin{equation}
  \dot x  = A(t) x
\label{eq:eqlin}
\end{equation}
where $x \in \mathbb{R}^m$ and $A:\R\to\R^{m\times m}$ is a bounded and continuous matrix function. Recall that solutions of~(\ref{eq:eqlin}) can be
written in terms of the evolution operator $\Phi:\R\times \R\to\R^{m\times m}$; the solution for the initial condition $x(t_0)=x_0$ is given by
\[
  t\mapsto \Phi(t,t_0) x_0\,.
\]
\begin{Definition}[uniform exponential stability]
  Consider the linear system~(\ref{eq:eqlin}) with evolution operator $\Phi$. System~(\ref{eq:eqlin}) is said to be \emph{uniformly exponentially stable} if there exists $K,\mu>0$ such that
  \begin{equation}
  \| \Phi(t,t_0) \| \le K e^{-\mu (t-t_0)}
  \quad \mbox{for all }\, t\ge t_0\,.
  \end{equation}
\end{Definition}

The following roughness theorem guarantees that uniform exponential stability is persistent under perturbations. A proof can be found in~\cite[Lecture~4, Prop.~1]{Coppel_78_1}.

\begin{Theorem}[roughness]\label{roughness}
  Consider the linear system~(\ref{eq:eqlin}) and assume that for $K>0$ and $\mu\in\R$, the evolution operator $\Phi$ satisfies the exponential estimate
  \begin{equation}\label{eq:exp-growth}
    \norm{\Phi(t,t_0)} \le K e^{- \mu(t-t_0)}
    \quad \mbox{for all }\, t\ge t_0\,.
  \end{equation}
  Consider a continuous matrix function $V:\R\to\R^{m\times m}$ such that
  \[
    \delta:=\sup_{t \in \R} \| V(t) \|<\infty\,.
  \]
  Then the evolution operator $\hat{\Phi}$ of the perturbed equation
  \[
    \dot  y = \big( A(t) + V(t) \big)y
  \]
  satisfies the exponential estimate
  \[
    \| \hat{\Phi}(t,t_0) \| \le K e^{- \hat{\mu}(t-t_0)}
    \quad \mbox{for all }\, t\ge t_0\,,
  \]
  where $\hat{\mu} := \mu - \delta K$.
\end{Theorem}

There are various criteria to obtain conditions for uniform exponential stability. We shall use the following criterion for diagonal dominant matrices, which can be found in~\cite[Lecture~6, Prop.~3]{Coppel_78_1}.

\begin{Propo}[diagonal dominance criterion]\label{prop:uniCont}
  Consider the linear system~(\ref{eq:eqlin}) with complex time-dependent coefficient matrices $A(t) = (A_{ij}(t))_{i,j=1,\dots,m}$, and suppose that there exists a constant $\mu >0$ such that
  \begin{equation}
    \Re(A_{ii}(t)) + \sum_{j=1,  \atop j\not= i}^m |A_{i j }(t)| \le - \mu < 0
    \quad \mbox{for all }\, t\in\R \mbox{ and } i\in\{1, \dots , m\}\,.
    \label{eq:Uii}
  \end{equation}
  Then the evolution operator $\Phi$ of~(\ref{eq:eqlin}) satisfies
  \[
    \| \Phi(t,t_0) \| \le K e^{-\mu(t-t_0)}
    \quad \mbox{for all }\, t\ge t_0\,.
  \]
 with $K=K(m)\ge 1$.
\end{Propo}

\section{Auxiliary results}\label{auxiliary}

In this section, we obtain various exponential estimates for orbits near the synchronisation manifold $S$ of~(\ref{eq:md1}).
First, we introduce a convenient splitting of coordinates along the synchronisation manifold and complementary to it, and derive the equations with respect to these coordinates.
Then we prove linear stability of the synchronisation manifold. Here we distinguish between diagonalisable and non-diagonalisable Laplacians. The latter case will follow from approximation results on diagonalisable Laplacians and roughness of the exponential estimates.
Finally, we introduce the concept of a tubular neighbourhood as a final ingredient to tackle the general proof of nonlinear stability.

In order to treat noncompact absorbing sets $U$ in Assumption~A1, we reformulate this assumption as follows.

\medskip

\noindent \textbf{Assumption A1'.} \textsl{
The function $f$ is continuous in the first argument and continuously differentiable in the second argument, and there exists
an open simply connected set $U \subset \mathbb{R}^m$ with $C^1$-boundary that is \emph{$\epsilon$-inflowing invariant} for some $\epsilon > 0$, i.e.~for all $x \in \partial U$ with inward-pointing normal vector $q_x$, we have
\begin{equation}\label{eq:eps-inflowing}
  \langle q_x , f(t,x) \rangle \ge \epsilon \fa t\in\R \mbox{ and } x\in \partial U\,.
\end{equation}
Moreover, there exists a $\Delta > 0$ such that the Jacobian $D_2 f$ is uniformly continuous and bounded on $B_\Delta(U):=\bigcup_{x\in U} \set{y\in\R^m:\|x-y\|<\Delta}$, i.e.~for some $\varrho>0$, we have
\begin{displaymath}
  \| D_2 f(t,x)\| \le \varrho \quad \mbox{for all } t \in \mathbb{R} \mbox{ and } x \in B_\Delta(U)\,.
\end{displaymath}}

\medskip

Note that if the closure $\bar{U}$ is compact, then uniformity of the inflowing invariance condition as well as the uniform continuity of
$D_2 f$ and existence of a bound $\varrho$ follow automatically. In the noncompact case, we require uniform bounds on the $\Delta$-enlarged neighbourhood $B_\Delta(U)$ for technical reasons.

We first obtain equations that govern the dynamics near the synchronisation manifold. Using a tensor representation, we can write the $n m$-dimensional system~(\ref{eq:md1}) equations by means of a single equation. To this end, define
\[
  X := \mbox{col}( x_1 , \dots , x_n )\,,
\]
where $\mbox{col}$ denotes the vectorisation formed by stacking the column vectors $x_i$ into a single
column vector. Similarly, define
\[
  F(t,X) := \mbox{col}( f(t,x_1) , \dots, f(t,x_n) )\,.
\]
We can analyse small perturbations away from the synchronisation manifold in terms of the tensor representation
\begin{equation}\label{eq:coord-split}
  X = \1 \otimes s + \xi\,,
\end{equation}
where $\otimes$ is the tensor product and $\1 = \mbox{col}(1,\dots,1)\in \mathbb{R}^n$, which is the eigenvector of $L$
corresponding to the eigenvalue zero. Note that $\1 \otimes s$ defines the diagonal manifold, and we view
$\xi$ as a perturbation to the synchronised state.

The state space $\R^n \otimes \R^m$ can be canonically identified with $\R^{nm}$, which we will use for shorter notation.
The coordinate splitting~(\ref{eq:coord-split}) is associated to a splitting of $\R^{nm}$ as the direct sum of subspaces
\begin{equation*}
  \R^{nm} = M \oplus N
\end{equation*}
with associated projections
\begin{equation*}
  \pi_M: \R^{nm} \to M, \qquad
  \pi_N:\R^{nm} \to N.
\end{equation*}
The subspaces $M,N \subset \R^{nm}$ are determined by embeddings from $\R^m$ and $\R^{(n-1)m}$, respectively, induced by the Laplacian
$L$ on $\R^n$.

Let us for the moment use the simplifying assumption that $L$ is diagonalisable with eigenvectors $\1, v_2, \dots, v_n$.
Then the subspaces $M,N$ have natural representations in terms of these eigenvectors as
\begin{equation*}
  M = \textrm{span}(\1) \otimes \R^m\,, \qquad
  N = \textrm{span}(v_2, \dots, v_n) \otimes \R^m\,.
\end{equation*}
This means that we have `natural' embeddings that induce coordinates on these subspaces:
\begin{eqnarray*}
  &\iota_M: \R^m \to M\,,
  &s \mapsto \1 \otimes s =  \mbox{col}(s,\ldots,s)\,,\\
  &\iota_N: \R^{(n-1)m} \to N\,,\quad
  &(y_2,\ldots,y_n) \mapsto \sum_{j=2}^n v_j \otimes y_j\,.
\end{eqnarray*}
If we drop the assumption that $L$ is diagonalisable, then we lose
the natural choice of an embedding for $N$. Note, however, that $N$ is
still determined as the eigenspace of all non-zero eigenvalues.

Note that the norm on $\R^{nm}$ we chose is the maximum over the Euclidean norm on $\R^m$, see~(\ref{eq:max-eucl-norm}). The norm $\norm{\cdot}$ on $\R^{nm}$ can be restricted to the subspaces $M,N$ and induces norms on the `coordinate' spaces $\R^m$ and $\R^{(n-1)m}$ by pullback under the embeddings.
Then the induced norm on $s \in \R^m$ is given by
\begin{equation}\label{eq:M-isometry}
  \norm{s}_{\iota_M} = \norm{\iota_M(s)} = \norm{\1 \otimes s}\,,
\end{equation}
which is precisely the Euclidean norm. Similarly, $\iota_M$ induces an
inner product on $M$. Henceforth, we shall identify $s \in \R^m$ with
$\1 \otimes s \in M$ under the isometry $\iota_M$.

Using the representation~(\ref{eq:coord-split}) for $X\in\R^{nm}$, given an initial condition $X_0=(s_0,\xi_0)$, the corresponding solution to~(\ref{eq:md1}) reads as $X(t)=(s(t),\xi(t))$. In the next result, we derive differential equations for these two components in a neighbourhood of the synchronisation manifold.

\begin{Propo}\label{prop:ODE-xi}
The two components of the solution $X(t)=(s(t),\xi(t))$ satisfy the system of equations
\begin{eqnarray}
  \1\otimes\dot{s} &=& \1\otimes f(t,s) + R_s(s,\xi)\,, \label{eq:ODE-s}\\
  \dot{\xi}        &=& T(t,s) \xi   + R_{\xi}(s,\xi)\,, \label{eq:ODE-xi}
\end{eqnarray}
where
\begin{equation}\label{eq:T-linear}
  T(t,s) = I_n \otimes D_2 f(t, s) - \alpha (L \otimes \Gamma)
\end{equation}
and $R_\ast:= R_s, R_\xi$ are the remainder functions such that for any $\epsilon > 0$, there is a $\delta > 0$ such that
for all $\| \xi \| \le \delta$, one has $\norm{R_\ast(s,\xi)} \le \epsilon \norm{\xi}$.
\end{Propo}

\begin{proof}
By Assumption~A2, Taylor's theorem implies that given $\epsilon>0$, there exists a $\delta>0$ such that
\[
  h(x) = \Gamma\,x + r(x) \quad  \mbox{with } \norm{r(x)} \le \epsilon \norm{x} \mbox{ whenever } \norm{x} \le \delta\,.
\]
Now we define
\begin{eqnarray*}
R_h(X)_i
&=& \sum_{j=1}^n W_{ij} r(x_i - x_j) = \sum_{j=1}^n W_{ij} r\big(p_i(\1 \otimes s + \xi) -
                              p_j(\1 \otimes s + \xi)\big)\\
&=& \sum_{j=1}^n W_{ij} r\big(p_i(\xi) - p_j(\xi)\big)\,,
\end{eqnarray*}
where $p_i: \R^{nm}\to\R^m$ maps canonically to the $i$-th component of the argument, $i\in\set{1,\dots,n}$. The vectors $R_h(X)_i \in\R^m$, $i\in\set{1,\dots, n}$ define a vector in $\R^{nm}$. Note that $R_h(X) = R_h(\xi)$ does not depend on $s \in M$ and satisfies the estimate
\[
  \norm{R_h(\xi)}
  \le \max_{i=1,\dots,n} \Bigg(\sum_{j=1}^n \abs{W_{ij}}\Bigg) \epsilon\,2 \norm{\xi}
  \qquad\textrm{whenever } \norm{\xi} \le \textstyle\frac{\delta}{2}\,.
\]
Recall that $L_{ij} = \delta_{ij} V_{i} - W_{ij}$, so the coupling term can then be rewritten as
\begin{equation}\label{EqR}
   \sum_{j=1}^{n} W_{ij} h(x_j - x_i)
= -\sum_{j=1}^{n} L_{ij} \Gamma\,x_j + R_h(\xi)_i
\end{equation}

The Taylor expansion of $F(t,X)$ around $\1 \otimes s$ reads as
\begin{eqnarray*}
  F(t,\1 \otimes s + \xi)
  &=& F(t,\1 \otimes s)
    + D_2 F(t,\1 \otimes s) \xi
    + R_F(t,s,\xi)\\
  &=& \1 \otimes f(t,s)
    + I_n \otimes D_2 f(t,s) \xi
    + R_F(t,s,\xi),
\end{eqnarray*}
where $\norm{R_F(t,s,\xi)} \le \epsilon\norm{\xi}$
when $\norm{\xi} \le \delta$. An algebraic manipulation of~(\ref{EqR}) allows a representation in coordinates $(s,\xi) \in M \oplus N$ of the $n$ equations forming~(\ref{eq:md1}):
\begin{eqnarray}
  \dot X = \1 \otimes \dot s + \dot \xi
  &=& \1 \otimes f(t,s)
    + I_n \otimes D_2 f(t,s) \xi
    - \alpha (L \otimes \Gamma) \xi \nonumber\\
 && + R_F(t,s,\xi)
    + \alpha R_h(\xi), \label{eq:ODE-X}
\end{eqnarray}
where we used $L\,\1 = 0$. Hence, the term $(L \otimes \Gamma)(\1 \otimes s)$ vanishes.

Next, we project the differential equation~(\ref{eq:ODE-X}) onto the spaces $M$ and $N$ to obtain
differential equations for $s$ and $\xi$:
\begin{eqnarray*}
  \1\otimes\dot{s} &=& \1\otimes f(t,s) + \pi_M(R_F(t,s,\xi) + \alpha R_h(\xi)),\\
  \dot{\xi}
      &=& T(t,s) \xi + \pi_N(R_F(t,s,\xi) + \alpha R_h(\xi)),
\end{eqnarray*}
where
\[
  T(t,s) = I_n \otimes D_2 f(t,s) -\alpha (L \otimes \Gamma).
\]
Note that both $I_n \otimes D_2 f(t,s)$ and
$L \otimes \Gamma$ preserve the subspaces $M$ and $N$, since
$I_n$ and $L$ preserve both $\textrm{span}(\1)$ and
$\textrm{span}(v_2, \ldots, v_n)$, so the projections can
be dropped there.
\end{proof}

\subsection{Diagonalisable Laplacians}\label{sec:diagL}

We now prove stability of the linear flow~(\ref{eq:T-linear}) for
$\xi \in N$, along any curve $s(t) \in S$, which is not necessarily a
solution. We first treat the diagonalisable case, and then the
non-diagonalisable one. Then, in Section~\ref{sec:sync}, we use
these results to prove stability of the fully nonlinear problem.

\begin{Lemma}[Diagonalisable case]
\label{lem:NetDiag}
Consider the linearisation of~(\ref{eq:ODE-xi}), given by
\begin{equation}\label{eq:linxi}
  \dot \xi = T(t, s(t))\xi\,, \quad \xi \in N
\end{equation}
with $s(t) \in U$, and the representations
\[
  L = P \Lambda P^{-1} \quad\mbox{and}\quad \Gamma = Q B Q^{-1}
\]
with $P\in\R^{n\times n}$ and $Q\in\R^{m\times m}$, such that $\Lambda = \mathrm{diag}(\lambda_1, \lambda_2, \dots, \lambda_n)$ and
$B = \mathrm{diag}(\beta_1,\dots,\beta_m)$. Then there exists a $\rho>0$ such that for all coupling strengths
\[
\alpha > \frac{\rho}{\gamma}\,,
\]
the evolution operator $\Phi$ of~(\ref{eq:linxi}) satisfies the estimate
\begin{equation*}
  \| \Phi(t,t_0) \| \le K \kappa(P \otimes Q)\,e^{-(\alpha\gamma - \rho)(t-t_0)} \fa t\ge t_0\,,
\end{equation*}
with $K \ge 1$, and where $\kappa(P \otimes Q)$ denotes the conditional number of $P \otimes Q$.
\end{Lemma}

Note that for matrices $P \in \mathbb{R}^{n \times n}$ and $Q \in \mathbb{R}^{m \times m}$,
we obtain
\begin{displaymath}
  \norm{P \otimes Q} \le  \norm{P}_{\infty} \norm{Q}_2\,,
\end{displaymath}
which implies that $\kappa(P \otimes Q) \le \kappa_{\infty}(P) \kappa_2(Q)$.

\begin{proof}[Proof of Lemma~\ref{lem:NetDiag}]
Note that $O := P \otimes Q$ is an invertible matrix that diagonalises $L \otimes \Gamma$, and the change of coordinates
\begin{equation}\label{eq:K-blockdiag}
  \tilde{T}(t)
  = O^{-1}\,T(t,s(t))\,O
  = I_n \otimes Q^{-1}\,D_2 f(t,s(t))\,Q
   -\alpha\,\Lambda \otimes B
\end{equation}
reduces $T(t)$ to $m$-block diagonal form. Thus, we have
\begin{equation*}
  \tilde{T} (t)= \bigoplus_{i=1}^n \tilde{T}_i(t) = \mathrm{diag} \big(\tilde T_1(t),\dots,\tilde T_n(t)\big)\,,
\end{equation*}
where
\[
\tilde{T}_i(t) := \underbrace{Q^{-1}\,D_2 f(t,s(t))\,Q}_{\tilde A(t):=} - \alpha\,\lambda_i\,B \fa t\in\R\,.
\]
Since for all $t\in\R$, the matrix $\tilde{T}(t)$ is block diagonal, the dynamics given by $\dot Y = \tilde T(t) Y$ preserves the splitting
$\R^{nm} = \bigoplus_{i=1}^n \R^m$, and hence, its associated evolution operator $\tilde{\Phi}$ is also of the form
\begin{equation}\label{eq:evolution-split}
  \tilde{\Phi}(t,t_0) = \bigoplus_{i=1}^n \tilde{\Phi}_i(t,t_0) \fa t,t_0\in\R\,,
\end{equation}
where each $\tilde{\Phi}_i$ is the evolution operator of $\dot y_i = \tilde T_i(t)y_i$. Note that restricting $T$ to $N$ corresponds to restricting $\tilde{T}$ to the blocks $i \ge 2$. The dynamics in each block is determined by
\begin{equation}
 \dot y_i =  ( \tilde{A}(t) - \alpha\,\lambda_i\,B) y_i\,.
\end{equation}
Now define
\begin{displaymath}
\tilde \rho :=  \sup_{t \in \R,\, s \in U} \big\| \tilde{A}(t) \big\|\,.
\end{displaymath}
Note that the matrix $\tilde A(t)$ depends implicitly on $s(t) \in U$, so by Assumption~A1 we get the estimate
\begin{equation}\label{rhovarrho}
  \tilde \rho \le \kappa(Q)\varrho\,.
\end{equation}
To apply Proposition~\ref{prop:uniCont}, we search for a condition on $\alpha$ such that
\begin{equation}\label{useofddc}
  \Re \big(\tilde{A}_{kk} - \alpha \lambda_i \beta_{k}\big)+\sum_{1 \le j \le m\atop j\not=k} \big|\tilde A_{kj}(t)\big| < 0 \fa k \in\set{1,\dots,m}\,.
\end{equation}
Since $\Re(\tilde{A}_{kk}) \le |\tilde{A}_{kk}|$, it is therefore sufficient that
\[
  \alpha > \frac{\sum_{j=1}^m |\tilde{A}_{kj}|}{\Re(\lambda_i \beta_{k})}
\]
holds. Note that $\Re(\lambda_i \beta_{k}) \ge \gamma$, so if we define
\begin{displaymath}
\sum_{j=1}^m |\tilde{A}_{ij}| \le c \tilde \rho =: \rho\,,
\end{displaymath}
where $c>0$ depends on the choice of the norm.
Then by the diagonal dominance criterion (Proposition~\ref{prop:uniCont}), the evolution operator $\tilde{\Phi}_i$ satisfies
\begin{equation}\label{eq:Tdiag}
  \| \tilde{\Phi}_i (t,t_0)\| \le K e^{-( \alpha \gamma - \rho ) (t-t_0)} \fa t \ge t_0\,.
\end{equation}
Finally, using~(\ref{eq:evolution-split}) and changing back to the
original coordinates, we have
\begin{eqnarray}
  \norm{\Phi(t,t_0)}
  &=& \big\|O\big(\textstyle\bigoplus_{i \ge 2} \tilde{\Phi}_i(t,t_0)\big)O^{-1}\big\| \nonumber\\
  &\le& \kappa(O)\,\textstyle \max_{i \ge 2} \big\|\tilde{\Phi}_i(t,t_0)\big\| \nonumber\\
  &\le& K \kappa(O)\,e^{-(\alpha \gamma - \rho)(t-t_0)}  \fa t \ge t_0\,.
  \label{eq:evolution-growth}
\end{eqnarray}
Note that $O^{-1}$ maps $M$ and $N$ onto the first and last $n-1$
of the $m$-tuples in $\R^{nm}$ respectively, so the restriction to
$N$ reduces to a direct sum over $i \ge 2$ after conjugation with
$O$, while we can simply estimate $\kappa(O|_{O^{-1}N}) \le \kappa(O)$.
\end{proof}

\subsection{Non-diagonalisable Laplacian}

We now treat the case when the Laplacian is non-diagonalisable and $\Gamma$ is diagonalisable. Note that if  $\Gamma$ is non-diagonalisable, the results follow from the density of diagonalisable matrices and the roughness property.

\begin{Lemma}[Non-diagonalisable Laplacian]
\label{lem:Nondiag}
Consider the situation of Lemma~\ref{lem:NetDiag} without the condition that the Laplacian is diagonalisable. Then there exists a $\bar\rho>0$ such that for all coupling strengths
\[
\alpha > \frac{\bar\rho}{\gamma}\,,
\]
the evolution operator $\Phi$ of~(\ref{eq:linxi}) satisfies the
estimate
\begin{equation*}
  \| \Phi(t,t_0) \| \le \bar C  e^{-(\alpha\gamma - \bar\rho)(t-t_0)} \fa t\ge t_0\,,
\end{equation*}
where $\bar C=\bar C(\Gamma, L) \ge 1$.
\end{Lemma}

The proof of this lemma makes use of roughness of exponential dichotomies and the density of diagonalisable Laplacians. We first establish the following auxiliary result.

\begin{Propo}\label{prop:HAp}
Let $\epsilon>0$ and $J$ be a complex Jordan block of dimension $m$. Consider
\[
\tilde J = J + E
\]
where $E = \mathrm{diag}(0,  \epsilon, 2 \epsilon, \dots ,  (m-1) \epsilon)$. Then there exists an $R \in \R^{m\times m}$ such that $R^{-1} \tilde J R$ is diagonal and
\[
\| R^{-1} E R \| = b \epsilon
\]
with the constant $b =b(m)$.
\end{Propo}

\begin{proof}
Note that $\tilde J$ is diagonalisable, since all the eigenvalues are distinct. All corresponding transformations $R$ are matrices of eigenvectors, upper triangular and can be computed explicitly. We normalise the eigenvectors such that for $\ell, j \in\set{1,\dots, m}$
\begin{displaymath}
  R_{\ell j} :=
 \left\{
 \begin{array}{c}
  \frac{(j-1)!}{(j-\ell)!} \epsilon^{\ell-1} \fa \ell \le  j \, , \\
  0 \mbox{ ~ otherwise }
\end{array}
\right.
\end{displaymath}
It is easy to verify that the elements $R^{-1}_{ik}$ with $i, k \in\set{1,\dots, m}$ of the inverse of $R$ read as
\begin{displaymath}
  R^{-1}_{ik} =
 \left\{
 \begin{array}{c}
    \frac{(-1)^{i+k}}{(i-1)!(k-i)!} \epsilon^{-(k-1)}  \fa i \le  k \, , \\
    0 \mbox{ ~  otherwise  }
\end{array}
\right.
\end{displaymath}
We have
\begin{displaymath}
(R^{-1} E R)_{ij} = \sum_{k,\ell} R^{-1}_{ik} E_{k\ell} R_{\ell j} \\
= \epsilon \frac{(-1)^{i}(j-1)!}{(i-1)!}\sum_{k=i}^{j} \frac{(-1)^{k}(k-1)}{(j-k)!(k-i)!}.
\end{displaymath}
Note that $(R^{-1} E R)_{ii} = (i-1)\epsilon$. Likewise, we have $(R^{-1} E R)_{i,i+1} = -  i \epsilon$.
Moreover, if $j > i+1$ then $(R^{-1} E R)_{ij} = 0$, since
\[
\sum_{k=i}^{j} \frac{(-1)^{k}(k-1)}{(j-k)!(k-i)!} = (-1)^i \sum_{l=1}^{j-i}\frac{(-1)^{l}}{(l-1)!(j-i-l)!} = 0\,.
\]
Therefore, $\max_{1\le i \ne m}  \sum_{j=1}^{m}| (R^{-1} E R)_{ij}| =\max\{ (2m-3),m-1 \}\epsilon$, and the result follows.
\end{proof}

Now we are ready to prove our approximation result.

\begin{Propo}\label{prop:GL}
Let $L$ be a Laplacian with simple eigenvalue zero and $\1$ its associated eigenvector. Then for any $\epsilon >0$, there exists
a matrix $\tilde L$ with simple eigenvalue zero and $\1$ its associated eigenvector such that
\begin{itemize}
\item[(i)] $\tilde L = P \tilde{\Lambda} P^{-1}$ with a diagonal matrix $\tilde{\Lambda} \in \R^{n\times n}$, and
\item[(ii)] $\big\| P^{-1} (\tilde L - L ) P \big\| \le \epsilon$.
\end{itemize}
\end{Propo}

\begin{proof}
We only need to prove the statement if $L$ is non-diagonalisable. We decompose $L$ in its complex Jordan canonical form
\[
  L = O J O^{-1}\,,
\]
where $J$ is a block diagonal matrix. The first block corresponds to the simple eigenvalue zero, so the first row contains only zeros, that is, $J =\mathrm{diag}(0,J_1,\dots,J_k)$, where $J_i$ are  Jordan blocks corresponding to
non-zero eigenvalues. Without loss of generality, we consider $k=1$.

Define $v := O^{-1} \1$. By hypothesis, we have $L \1 = 0$,  so
\begin{equation}
J v  = 0\,.  \label{Tn}
\end{equation}

As each Jordan block has its own invariant subspace, (\ref{Tn}) implies
$v = (1,0,\dots,0)$. Define $E := \mathrm{diag}(0,\epsilon, 2 \epsilon, \dots, (n-1) \epsilon)$, and note that
\begin{equation}
E v  = 0\,. \label{eq:lambda}
\end{equation}

Consider the matrix
\begin{displaymath}
  \tilde L = O (J + E) O^{-1}\,,
\end{displaymath}
which is diagonalisable. Moreover, by~(\ref{Tn}) and~(\ref{eq:lambda}), we obtain that
$\tilde L$ has zero as a simple eigenvalue with associated
eigenvector $\1$. By Proposition~\ref{prop:HAp}, we obtain
\begin{displaymath}
  J + E = R \tilde{\Lambda} R^{-1}\,,
\end{displaymath}
and hence the matrix $P = OR$ diagonalises $\tilde L$. For this reason,
\begin{displaymath}
P^{-1} (\tilde L - L) P = P^{-1} (O E O^{-1}) P = R^{-1} E R\,,
\end{displaymath}
and the result follows by Proposition~\ref{prop:HAp}.
\end{proof}

\begin{proof}[Proof of Lemma~\ref{lem:Nondiag}]
As in the diagonalisable case, we consider the linearised equation~(\ref{eq:linxi}) for $\xi \in N$ along any curve $s(t) \in U$. By Proposition~\ref{prop:GL}, there is a diagonalisable matrix $\tilde L$ in an arbitrary neighbourhood of the Laplacian $L$.
We rewrite~(\ref{eq:linxi}) as
\begin{equation}\label{eq:linxi-pert-nondiag}
  \dot{\xi} = \big[I_n \otimes D_2 f(t,s(t)) - \alpha \tilde{L} \otimes \Gamma\big] \xi + \alpha\big[(\tilde{L} - L)\otimes \Gamma\big] \xi \,.
\end{equation}
Note that this is a small perturbation of the same equation
with diagonalisable Laplacian $\tilde{L}$, so we can apply the results from Subsection~\ref{sec:diagL}.
Recall that $\Gamma = Q B Q^{-1}$ and $\tilde{L} = P \tilde{\Lambda} P^{-1}$ (see Proposition~\ref{prop:GL}).
Moreover, consider the change of variables $\zeta = (P^{-1} \otimes Q^{-1}) \xi$.
We obtain
\begin{equation}\label{eq:zetaT}
  \dot{\zeta} = \big[I_n \otimes Q^{-1} D_2 f(t,s(t)) Q - \alpha \tilde{\Lambda} \otimes B \big] \zeta + \alpha\big[P^{-1}(\tilde{L} - L)P\otimes B \big] \zeta \,.
\end{equation}
We treat $\alpha \big(P^{-1} (\tilde L - L) P \otimes B \big) \zeta$ as a perturbation of the equation
\begin{equation}\label{eq:zeta}
  \dot \zeta
  = \big(I_n \otimes Q^{-1} D_2 f(t,s(t)) Q - \alpha(\tilde{\Lambda} \otimes B)\big)
  \zeta\,.
\end{equation}
It follows from the proof of Lemma~\ref{lem:NetDiag} (see~(\ref{eq:Tdiag}) for details) that  the evolution operator $\tilde{\Phi}$ of (\ref{eq:zeta})  satisfies
\[
  \| \tilde{\Phi}(t,t_0) \| \le K e^{- (\alpha \gamma - \rho) (t-t_0)}\,,
\]
where $K$ does not depend on $n$ as (\ref{eq:zeta}) is block diagonal. Theorem~\ref{roughness} (the roughness theorem) implies that the condition
\begin{equation}\label{eq:rho-nondiag-est}
  \alpha \| P^{-1} (\tilde L - L) P \otimes B \| < \frac{\alpha \gamma - \rho}{K}
\end{equation}
leads to an exponential stability estimate for the perturbed equation~(\ref{eq:linxi-pert-nondiag}).
By Proposition~\ref{prop:GL} (ii), we can choose $\tilde L$ such that
$\| P^{-1} (\tilde L - L) P \| \le \epsilon / \| B\|$, so~(\ref{eq:rho-nondiag-est}) is satisfied if
taking $\epsilon < (\alpha \gamma - \rho)/ (\alpha K)$.
Hence, setting $\bar{\rho} := \rho + \alpha K  \epsilon$, then for all $\alpha > \bar{\rho}/\gamma$ the linear flow $\Phi(t,t_0)$ for~(\ref{eq:linxi-pert-nondiag}) satisfies
\[
\| \Phi(t,t_0) \| \le K \kappa(P \otimes Q) e^{- (\alpha \gamma - \bar{\rho})(t-t_0)} \fa t\ge t_0\,,
\]
where the conditional number is due to transforming back to the original variables $\xi$.
\end{proof}

To analyse the solution curves $(s(t),\xi(t))$ of the nonlinear system~(\ref{eq:ODE-s},\ref{eq:ODE-xi}) we introduce the concept of a tubular neighbourhood.

\begin{figure}[htbp] 
   \centering
   \includegraphics[width=2.7in]{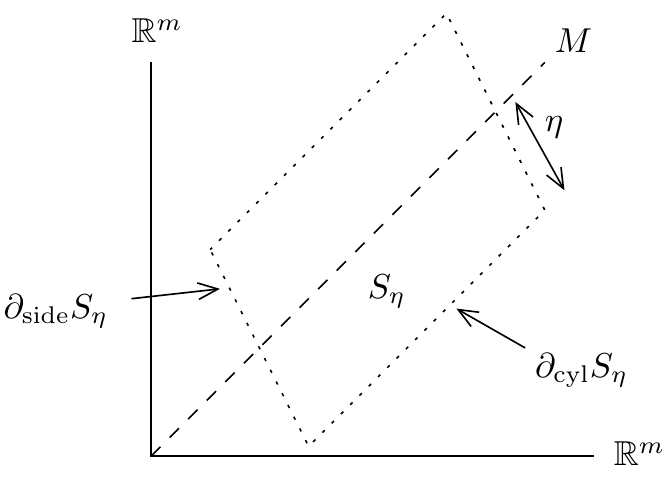}
  \caption{Tubular neighbourhood $S_\eta$ for $n = 2$.}
  \label{fig:tubular-neigh}
\end{figure}

\begin{Definition}[$\eta$-tubular neighbourhood]\label{def:tubular-neigh}
  Let $S = \1 \otimes U \subset M$ be a subset of the diagonal manifold. Then the set
  \begin{equation}\label{eq:tubular-neigh}
    S_\eta = \big\{ \1 \otimes s + \xi : s \in U \mbox{ and } \xi \in N, \mbox{ where }\norm{\xi} < \eta \big\}
  \end{equation}
  for a given $\eta>0$ is called the \emph{$\eta$-tubular neighbourhood} of $S$.
\end{Definition}
See Figure~\ref{fig:tubular-neigh} for a schematic illustration of this definition. Note that the directions along $N$ in which the tubular stretches out do not need to be orthogonal to $M$.

%

Assumption~A1' says that the single-node system has a uniformly inflowing invariant set $U \subset \R^m$. A similar result holds in a neighbourhood of the synchronisation manifold $S$ in the coupled network, since the following lemma implies that if the solution curve $(s(t),\xi(t))$ leaves $S_\eta$, then it must do so by $\norm{\xi(t)}$ growing larger than $\eta$.

\begin{Lemma}\label{lem:inflow-mfld}
  Consider  Assumption~A1' with the $\epsilon$-inflowing invariant set $U \subset \R^m$. Let
  $\dot{X} = F(t,X)$ describe the dynamics of $n$ uncoupled copies of this system and let
  $G: \R \times \R^{nm} \to \R^{nm}$ be a perturbation such that for some $r > 0$ and $\delta > 0$, one has
  \begin{equation*}
    \sup_{t \in \R, X \in S_r} \norm{G(t,X)}
    \le \delta < \frac{\epsilon}{\norm{\pi_M}}\,.
  \end{equation*}
  Then there exists an $\eta \in (0,r]$ such that solution curves
  $(s(t),\xi(t))$ of $\dot X = F(t,X) + G(t,X)$ can only leave the tubular neighbourhood $S_\eta$ through
  \begin{equation*}
    \partial_{\rm cyl} S_\eta
    := \{ \1 \otimes s + \xi :
         s \in U \mbox{ and } \norm{\xi} = \eta \}\,.
  \end{equation*}
\end{Lemma}

\begin{proof}
  Choose $\eta$ such that $0 < \eta \le r$. The boundary of $S_\eta$ consists of two parts:
\begin{equation*}
  \partial S_\eta = \partial_{\rm cyl} S_\eta \cup \partial_{\rm side} S_\eta\,,
\end{equation*}
where $\partial_{\rm side} S_\eta := \set{\1 \otimes s + \xi: \norm{\xi}\le\eta\mbox{ and } s \in \partial U}$.

We consider the dynamics on $\partial_{\rm side} S_\eta$. Let $q$ be
the inward pointing normal vector at $s \in \partial U$. Locally we
have $\partial_{\rm side} S_\eta = \partial S \oplus N$, so $F+G$
points inwards at $\1 \otimes q + \xi$ precisely if its projection
onto $M$ along $N$ has positive inner product with $q$. Note that we use the
isometry $\iota_M$ from~(\ref{eq:M-isometry}) to endow $M$ with the
inner product $\langle \,\cdot\,,\,\cdot\,\rangle_M$ induced from
$\langle \,\cdot\,,\,\cdot\,\rangle_{\R^m}$, but no inner product on
$\R^{nm}$ is used (nor defined).

If $\eta$ is chosen sufficiently small, then $S_\eta$ is contained
within the product space $B_\Delta(U)^n$ where we have uniform bounds
$\norm{D_2 F} \le \varrho$ and $\norm{G} \le \delta$. It follows that
\begin{eqnarray*}
  &&  \langle \iota_M(q),\pi_M[ F(t,\1 \otimes s + \xi)
                               +G(t,\1 \otimes s + \xi)]\rangle_M\\
  &=& \langle q,f(t,s)\rangle_{\R^m}
     +\langle q,\iota_M^{-1}\circ\pi_M[D_2 F(t,\1 \otimes s + \tau\,\xi)\xi
                                         + G(t,\1 \otimes s + \xi)]\rangle_{\R^m}\\
  &\ge& \epsilon - \norm{\pi_M}(\norm{D_2 F} \norm{\xi} + \norm{G})\\
  &\ge& \epsilon - \norm{\pi_M}(\varrho \eta + \delta),
\end{eqnarray*}
where we applied the mean value theorem with $\tau \in (0,1)$ as
interpolation variable. Since $\norm{\pi_M}\delta < \epsilon$,
there exists an $\eta > 0$ such that $F+G$
points inwards everywhere at $\partial_{\rm side} S_\eta$.
\end{proof}

Finally, we shall make use of the following lemma, which is a variant
on Gronwall's Lemma.
\begin{Lemma}\label{lem:varconst-est}
  Let $x(t) \in \R$ satisfy the integral inequality
  \begin{equation}\label{eq:varconst-ineq}
    x(t) \le C e^{-\mu(t-t_0)} x_0
            +\int_{t_0}^t C e^{-\mu(t-\tau)}\big(\alpha x(\tau) + \beta\big) \d\tau\,,
  \end{equation}
  with $C,\mu > 0$ and $x_0,\alpha,\beta \ge 0$, whenever $x \le \delta$.

  If $\tilde{\mu} := \mu - C \alpha > 0$ and
  $x_0 < \frac{1}{C}\big(\delta - \frac{\beta}{\tilde{\mu}}\big)$, then
  $x(t)$ is bounded by
  \begin{equation}\label{eq:varconst-sol}
    x(t) \le C e^{-\tilde{\mu}(t-t_0)} \Big(x_0  - \frac{\beta}{\tilde{\mu}}\Big)
            +\frac{C \beta}{\tilde{\mu}}  \fa t \ge t_0\,,
  \end{equation}
  and in particular $x(t) < \delta$ holds for all $t \ge t_0$.
\end{Lemma}

\begin{proof}
  The integral inequality is equivalent to the differential inequality
  \begin{equation*}
    \dot{x}(t) \le -\mu x(t) + C\big(\alpha x(t) + \beta\big)\,, \qquad x(t_0) = C x_0\,,
  \end{equation*}
  so by a standard application of Gronwall's lemma we
  obtain~(\ref{eq:varconst-sol}), as long as the solution satisfies
  $x(t) \le \delta$. Now assume by contradiction that this assumption
  is violated. Then there exists a $t_1 \ge t_0$ such that
  $x(t) = \delta$ for the first time at $t = t_1$. However, the
  assumption $x(t) \le \delta$ is true up to time $t_1$, so by the
  previous estimates and the assumption that
  $x_0 < \frac{1}{C}\big(\delta - \frac{\beta}{\tilde{\mu}}\big)$ it
  follows that $x(t_1) < \delta$. This contradiction completes the
  proof.
\end{proof}

\section{Synchronisation}\label{sec:sync}

In the previous section we have established all auxiliary results to prove our main theorem on synchronisation (Theorem~\ref{MainThm}), which will be restated for convenience.

\begin{theorem}[synchronisation]
Consider the network of diffusively coupled equations~(\ref{eq:md1}) satisfying~A1--A3. Then there exists a $\rho = \rho(f,\Gamma)$ such that for all coupling strengths
\[
  \alpha > \frac{\rho}{\gamma} \,,
\]
the network is \emph{locally uniformly synchronised}. This means that there exist a $\delta>0$ and a $C=C(L,\Gamma)>0$ such that if $x_i(t_0) \in U$ and $\| x_i (t_0) - x_j (t_0)\| \le \delta$ for any $i,j\in\{1,\dots,n\}$, then
\begin{displaymath}
  \| x_i (t) - x_j (t)\| \le C e^{-(\alpha \gamma - \rho)(t-t_0)} \| x_i (t_0) - x_j (t_0)\| \quad \mbox{for all } t\ge t_0\,.
\end{displaymath}
\end{theorem}

\begin{proof}
  Set
  \begin{equation*}
    X(t_0) = \1 \otimes s(t_0) + \xi(t_0)
          := \mbox{col}\big(x_1(t_0),\cdots,x_n(t_0)\big)
  \end{equation*}
  where $x_i(t_0) \in U$ and $U \subset \R^m$ is
  $\varepsilon$-inflowing invariant. Due to the uniformity assumptions
  in A1', there exists a slightly enlarged neighbourhood
  $B_{\Delta/2}(U) \supset U$ that is still $\varepsilon/2$-inflowing
  invariant. We set $S = \1 \otimes B_{\Delta/2}(U)$. If we choose the
  distance bound $\norm{x_i(t_0) - x_j(t_0)} \le \delta$ sufficiently
  small (depending on the angle between $M$ and $N$), then
  $s(t_0) \in S$ holds, while we also have
  $\norm{\xi(t_0)} \le \norm{\pi_N}\delta$.

  By Lemma~\ref{lem:inflow-mfld} there exists a tubular neighbourhood
  $S_\eta$ of positive size $\eta > 0$ over $S$ that is inflowing
  invariant on the `side' and contained within
  $B_\Delta(U)^n \subset \R^{nm}$, so the uniform assumptions of A1'
  hold.

  Now lemmas~\ref{lem:NetDiag} and~\ref{lem:Nondiag} together imply
  that there exists a $\rho > 0$ such that for
  $\alpha > \frac{\rho}{\gamma}$, the evolution operator $\Phi(t,t_0)$
  for $\xi$ satisfies an exponential estimate with decay rate
  $-(\alpha \gamma - \rho)$. The nonlinear remainder of the flow of
  $\xi$ can be bounded by an arbitrarily small linear term when
  $\norm{\xi}$ is small, as controlled by $\eta$. By variation of
  constants, Eq.~(\ref{eq:ODE-xi}) for $\xi$ is equivalent to
  \begin{equation}\label{eq:var-const-xi}
    \xi(t) = \Phi(t,t_0)\xi(t_0) + \int_{t_0}^t \Phi(t,\tau) R_\xi(s(\tau),\xi(\tau)) \d\tau\,.
  \end{equation}
  Now we assume that $\norm{\xi(t)} \le \eta$ for all $t \ge t_0$ and
  estimate
  \begin{equation*}
    \norm{\xi(t)}
    \le C e^{-(\alpha\gamma - \rho)(t-t_0)} \norm{\pi_N} \delta
         +\int_{t_0}^t C e^{-(\alpha\gamma - \rho)(t-\tau)} \epsilon(\eta) \d\tau\,.
  \end{equation*}
  Hence, when we choose $\delta < \frac{\eta}{C \norm{\pi_N}}$ and
  $\epsilon(\eta)$ sufficiently small, then we can apply
  Lemma~\ref{lem:varconst-est} with $\beta = 0$ and conclude that
  \begin{equation*}
    \norm{\xi(t)} \le C e^{-\tilde{\mu}(t-t_0)} \norm{\pi_N} \delta
    \fa t \ge t_0\,,
  \end{equation*}
  with $\tilde{\mu} = \alpha\gamma - \rho - C \epsilon(\eta)$. Thus,
  if we choose $\tilde{\rho} = \rho + C \epsilon(\eta)$, then for all
  $\alpha > \frac{\tilde{\rho}}{\gamma}$ the complete solution curve
  $(s(t),\xi(t))$ for the nonlinear system is contained in $S_\eta$
  for all $t \ge t_0$ and converges to the synchronisation manifold
  $S$ with decay rate $-(\alpha \gamma - \tilde{\rho})$. The explicit
  estimate for $\norm{x_i(t) - x_j(t)}$ can be recovered from
  \begin{equation*}
          \norm{x_i(t) - x_j(t)}
    \le 2 \norm{x_i(t) - s(t)}
    \le 2 \norm{\xi(t)}
  \end{equation*}
  and the fact that $\delta$ can be chosen smaller to match
  $\norm{x_i(t) - x_j(t)}$.
\end{proof}

\begin{Remark}\label{rem:delta-est}
  Explicit estimates for the size of $\delta$ in Theorem~\ref{MainThm}
  can be found when more details of the system are known. For example,
  if the second derivative of $f$ is bounded, i.e.
  \[
    \norm{D_2^2 f(t,x)} \le \sigma \fa t \in \mathbb{R} \mbox{ and } x \in U\,,
  \]
  and the coupling function is linear, i.e.~$h(x) = \Gamma x$, then
  $\delta$ can be estimated as
  \begin{equation}\label{eq:delta-est}
    \delta = \frac{\alpha\gamma - \rho}{4 \sigma C \norm{\pi_N}}\,.
  \end{equation}
  Note that for convenience, we ignore effects on the size of $\delta$
  introduced by estimates at the boundary of the synchronisation manifold.
  Under these assumptions the remainder $R_\xi$
  in~(\ref{eq:var-const-xi}) consists of $R_F$, the nonlinearities of $f$, and can be estimated as
  $\norm{R_F(t,s,\xi)} \le \sigma \norm{\xi}^2$ using mean value
  theorem arguments. To conclude the argument, fix $\delta = \eta/(2 C \norm{\pi_N})$ and
  follow the proof of Theorem~\ref{MainThm}.
\end{Remark}

\subsection{Behaviour of $\rho$ as function of $\Gamma$}\label{Aprho}

Our approach is constructive and allows to estimate the bounds for $\rho = \rho(f,\Gamma)$ whenever specific information on the function $h$ is provided. By Lemma~\ref{lem:Nondiag}, it is clear that the diagonalisation properties of the Laplacian have no effect on the bounds for $\rho$. In the following, we only discuss symmetric Laplacians $L$. As an illustration, we look at two cases for $\Gamma$.

\begin{enumerate}
\item[(i)] \textbf{$\Gamma$ is symmetric.} There exists an orthogonal matrix $Q$ such that $\Gamma = Q B Q^{-1}$. Note that $\kappa(Q) = 1$ (i.e.~the conditional number with respect to the Euclidean norm). From~(\ref{rhovarrho}), it follows that
\[
  \rho \le \tilde c \varrho
\]
for some $\tilde c>0$. The bound for $\rho$ is independent of $\Gamma$ for this reason. Note that this can be observed in the left panel of Figure~\ref{fig:Results_Lorenz}.

\item[(ii)] \textbf{$\Gamma$ is non-diagonalisable.} To treat the non-diagonalisable case, we employ the above perturbation techniques we developed for the Laplacian, i.e.~we approximate $\Gamma$ by a diagonalisable matrix $\tilde \Gamma$. Notice that $\Gamma$ can be represented in its Jordan form $\Gamma = Q J Q^{-1}$, and we can write $\tilde J = J + E$, where $E$ is an $\epsilon$-perturbation diagonal matrix as in Proposition~\ref{prop:HAp}. The approximation $\tilde \Gamma$ reads as $\tilde \Gamma = Q(J+E)Q^{-1}$, and as in Proposition~\ref{prop:HAp}, if $P$ denotes the matrix that diagonalises $J+E$ (i.e.~$\tilde B = P^{-1}(J+E)P$ is diagonal), then $\tilde \Gamma = QP \tilde B P^{-1}Q^{-1}$. Hence,
\begin{displaymath}
  \rho \le c \varrho \kappa(QP) \le c \varrho \kappa(Q)\kappa(P)\,.
\end{displaymath}
By Proposition~\ref{prop:HAp}, it is easy to check that
\begin{displaymath}
  \kappa(P) = \| P \|  \| P^{-1} \| \le   \frac{d}{\epsilon^{m-1}}\,,
\end{displaymath}
where $d>0$ does not depend on $\epsilon$. The aim is to minimise $\rho$, which means minimising $\kappa(P)$. The perturbation size $\epsilon$ should be of the same order as $\beta$, since the real parts of the eigenvalues of $J+E$ must be positive. This can be obtained, for instance, by choosing $\epsilon = r \beta$ for some fixed $r\in(0,1)$. This yields to the following bound
\[
\rho \le \frac{k}{\beta^{m-1}}\,,
\]
where $k$ is a constant.

\end{enumerate}

Note the different behaviour for the bound as a function of $\beta$ between the case when $\Gamma$ is symmetric and when $\Gamma$ is non-diagonalisable. This helps to explain the nonlinear behaviour observed in Figure~\ref{fig:LambdaAlpha} and in the right panel of Figure~\ref{fig:Results_Lorenz}.

\section{Persistence}\label{sec_pers}

As in the previous section, we make use of the auxiliary results from Section~\ref{auxiliary} in order to prove our main theorem on persistence (Theorem~\ref{PertThm}), which will be restated for convenience.

\begin{theorem}[persistence]
  Consider the perturbed network~(\ref{eq:perturbed}) of diffusively coupled equations fulfilling Assumptions A1--A3, and suppose that
  \[
    \alpha > \frac{\rho}{\gamma}
  \]
  as in Theorem~\ref{MainThm}. Then there exist $\delta>0$, $C>0$ and $\epsilon_g>0$ such that for all $\epsilon_0$-perturbations satisfying
  \begin{displaymath}
    \norm{g_i(t, x)} \le \epsilon_0 \le \epsilon_g \quad \mbox{for all } t\in\R\,,\, x\in U\mbox{ and } i\in\{1,\dots, n\}
  \end{displaymath}
  and initial conditions satisfying $\norm{x_i(t_0) - x_j(t_0)} \le \delta$ for any $i,j\in\{1,\dots,n\}$, the estimate
  \begin{displaymath}
    \hspace{-1cm}\norm{x_i(t) - x_j(t)}
    \le C e^{-(\alpha \gamma - \rho )(t-t_0)} \norm{x_i(t_0) - x_j(t_0)}
       +\frac{C\epsilon_0}{\alpha \gamma  - \rho}
    \quad \mbox{for all } t \ge t_0
  \end{displaymath}
  holds.
\end{theorem}

Note that the proof of this theorem does not specifically depend on the fact that the perturbations $g_i$ of the nodes are decoupled; the function $G$ below can depend arbitrarily on the total state $X$ (or can be subjected to random perturbations).

\begin{proof}[Proof of Theorem~\ref{PertThm}]
  Denote by
  \[
  G(t,X) = \mbox{col}( g_1(t,x_1) , \dots , g_n(t,x_n) )
  \]
  the perturbation for the network and note that $\norm{G} \le \epsilon_0$.
  As in the proof of Theorem~\ref{MainThm},
  Lemma~\ref{lem:inflow-mfld} guarantees that there exists an
  $\eta$-tubular neighbourhood $S_\eta$ such that solutions of the
  complete system for $(s,\xi)$ cannot escape along $s$,
  when $\epsilon_g,\eta$ are sufficiently small.

The perturbed network equation for $X=(s,\xi)$ in $S_\eta$ now reads as
\begin{displaymath}
  \dot X = F(t,X) - \alpha L\otimes \Gamma \xi + R_h(\xi)+G(t,X)\,,
\end{displaymath}
where $R_h$ is the Taylor remainder associated with the coupling function $h$. Projecting this equation onto the synchronisation manifold yields an equation for the component $s$ of $X$. On the other hand, the differential equation for $\xi$ is given by
  \begin{equation}\label{eq:ODE-xi-mod}
    \dot{\xi}
    = T(t,s(t)) \xi + R(t, s(t),\xi) + \pi_N(G(t, \1 \otimes s + \xi))\,,
  \end{equation}
  see Proposition~\ref{prop:ODE-xi}. Let $\epsilon(\eta)$ denote a Lipschitz constant within $S_\eta$ of $R$ with respect to $\xi$, which does not depend on $t$.

  In the same way as in the proof of Theorem~\ref{PertThm}, we obtain
  a variation of constants formula for solutions of~(\ref{eq:ODE-xi-mod}),
  \[
    \xi(t) = \Phi(t,t_0)\xi(t_0) + \int_{t_0}^t \Phi(t,\tau)
              \big[R(\tau,s(\tau),\xi(\tau))+\pi_{N}( G(\tau,\1 \otimes s(\tau) + \xi(\tau)))\big] \d\tau\,.
  \]
  With initial conditions $\norm{\xi(t_0)} \le \norm{\pi_N}\delta$,
  lemmas~\ref{lem:NetDiag} and~\ref{lem:Nondiag}, and the assumption
  that
  \begin{equation*}
    \norm{\xi(t)} \le \delta_1 < \eta \fa t \ge t_0\,,
  \end{equation*}
  this leads to the estimate
  \begin{equation*}
    \norm{\xi(t)}
    \le C e^{-\mu t} \norm{\pi_N}\delta
       +\int_{t_0}^t C e^{-\mu(t-\tau)} (\epsilon(\delta_1)\norm{\xi(\tau)} + \norm{\pi_N}\epsilon_0) \d\tau\,,
  \end{equation*}
  where $\mu = \alpha \gamma - \rho$. We choose
  $\delta < \frac{\eta}{C \norm{\pi_N}}$ and $\delta_1,\epsilon_g$
  sufficiently small and apply Lemma~\ref{lem:varconst-est} with
  $\alpha = \epsilon(\delta_1), \beta = \norm{\pi_N}\epsilon_0$ to
  find that
  \begin{equation}\label{eq:persist-convergence}
    \norm{\xi(t)}
    \le C e^{\tilde{\mu}(t-t_0)}\norm{\pi_N}
        \Big(\delta - \frac{\epsilon_0}{C\tilde{\mu}}\Big)
       +\frac{C \norm{\pi_N}\epsilon_0}{\tilde{\mu}} \fa t \ge t_0\,,
  \end{equation}
  where $\tilde{\mu} = \alpha \gamma - \rho - C \epsilon(\delta_1)$.
  As in the proof of Theorem~\ref{PertThm}, we choose
  $\tilde{\rho} = \rho + C \epsilon(\delta_1)$ instead of $\rho$ and
  the estimate for $\norm{x_i(t) - x_j(t)}$ follows
  from~(\ref{eq:persist-convergence}) by adapting $\delta$.
\end{proof}

In particular, note that asymptotically, the bound in~(\ref{eq:persist-convergence}) converges to $\frac{C \norm{\pi_N}\epsilon_0}{\alpha \gamma - \tilde{\rho}}$. Furthermore, it follows from the details of Lemma~\ref{lem:NetDiag} that the constant $C$ depends on the Laplacian $L$ only through its conditional number $\kappa(P)$.

Finally, we can proof Corollary~\ref{corofpersistence} from the Introduction.

\begin{proof}[Proof of Corollary~\ref{corofpersistence}]
This corollary is a direct consequence of our persistence result.
For simplicity, we now endow the space $\R^{nm}$ with the Euclidean norm
\[
\| X \|_2 = \Big( \sum_{i=1}^n \| x_i \|_{2}^2 \Big)^{1/2} \fa X = \mathrm{col}(x_1, \dots, x_n) \in \mathbb{R}^{nm}\,.
\]
Note that in view of~(\ref{eq:persist-convergence}), for large times, we obtain
\begin{equation}\label{eq:x}
\| \xi \|_2 = \left(  \sum_{i=1}^n \|s - x_i \|^2_2 \right)^{1/2} \le \frac{2 K \kappa_2(P \otimes Q) \| G\|_2}{\mu}
\end{equation}
where the contraction rate $\mu$ is given by $\mu = \alpha \gamma - \rho$. For simplicity, we omit the arguments of the functions $s$, $x$, $G$ and $\xi$.

Moreover, $\kappa_2(P\otimes Q) \le \kappa_2(P) \kappa_2(Q)$, and
since the Laplacian is symmetric, it can be diagonalised
by an orthogonal similarity transformation, which implies that $\kappa_2(P)=1$ together with $\|\pi_N\|_2=1$.
Moreover, by the equivalence of norms we obtain
\[
\| G \|_2 \le \sqrt{n} \| G \| \le \sqrt{n} \epsilon_0,
\]
Replacing this estimate in (\ref{eq:x}) we obtain
\begin{equation}\label{eq:x1}
\left(  \sum_{i=1}^n \|s - x_i \|^2_2 \right)^{1/2}  \le \frac{\tilde K \sqrt{n}\epsilon_0}{\mu} \,,
\end{equation}
where $\tilde K = 2 K \kappa_2(Q)$.  We scale equation~(\ref{eq:x1}) to obtain
\begin{equation}\label{eq:aux}
  \left( \frac 1 n  \sum_{i=1}^n \| s - x_i \|^2_2 \right)^{1/2} \le   \frac{\tilde K \epsilon_0}{\mu }\,,
\end{equation}
and applying the sum of squares inequality
\[
\frac{1}{n}\sum_{i=1} a_i \le \sqrt{\frac{1}{n} \sum_{i=1}^n a_i^2 }
\]
leads to
\begin{equation}\label{smean}
\frac 1 n \sum_{i=1}^n \| s - x_i \|_2  \le \frac{\tilde K \epsilon_0}{\mu }\,.
\end{equation}
The triangle inequality implies
\begin{eqnarray*}
\frac 1 n \left| \sum_{i=1}^n \| s -x_j \|_2 - \| x_j - x_i \|_2   \right| & \le & \frac 1 n \sum_{i=1}^n \left|  \| s -x_j \|_2 - \| x_j - x_i \|_2  \right| \\ &\le &   \frac 1 n \sum_{i=1}^n \|s - x_i \|_2\,.
\end{eqnarray*}
Hence,
\begin{displaymath}
\frac 1 n \left| \sum_{i=1}^n \left( \| s -x_j \|_2 - \| x_j - x_i \|_2  \right) \right| \le  \frac{\tilde K \epsilon_0}{\mu}\, ,
\end{displaymath}
as we control the first sum by (\ref{smean}) we obtain
\begin{displaymath}
  \frac{1}{n} \sum_{i=1}^n  \| x_j - x_i \|_2  \le  \frac{2 \tilde K \epsilon_0}{\mu}\,.
\end{displaymath}
To conclude the result, we take the sum over the index $j$ and divide by the network size $n$.
This finishes the proof of this corollary.
\end{proof}

\section{Generalisations}\label{secgeneral}

Although our set-up is very general and includes non-autonomous systems and non-diagonalisable Laplacians, the assumptions we make are only sufficient for synchronisation, but not necessary. For instance, let $(u,v) = x \in  \mathbb{R}^2$ and consider as isolated dynamics $\dot{x} = f(x)$ with $f(x) = (u,u-v)$, and
\[
\begin{array}{cc}
\dot{x}_1 = &  f(x_1) + \alpha \Gamma (x_2  - x_1) \\
\dot{x}_2 = &  f(x_2) + \alpha \Gamma (x_1  - x_2) \\
\end{array}
\quad\mbox{with}\quad
\Gamma =
\left(
\begin{array}{cc}
1 & 0 \\
0 &0 \\
\end{array}
\right).
\]
Note that in this situation $\Gamma$ has an eigenvalue zero, so Assumption~A3 is violated. However, this coupled system synchronises for $\alpha > 1/2$. This happens as all instabilities occurs due to the first variable, and the coupling $\Gamma$ acts solely on this variable. For a numerical example of a chaotic system displaying synchronisation with only one variable coupled, see~\cite{Pecora_98_1}.

The boundedness of the Jacobian $D_2 f$ in Assumption~A1', and Assumption~A3 are used in Lemma~\ref{lem:NetDiag} to obtain uniform exponential stability of the linear system (\ref{eq:linxi}). For this purpose, we use the diagonal dominance criterion, see (\ref{useofddc}) in the proof of Lemma~\ref{lem:NetDiag}. It is clear that one could get uniform exponential stability without the two above mentioned assumptions. Note that under reasonable assumptions, a necessary and sufficient condition for uniform exponential stability (and thus persistent synchronisation) is that the dichotomy spectrum of (\ref{eq:linxi}) is contained in the negative half line~\cite{Kloeden_11_2} (see~\cite{Froyland_13_1} for a comparative study of numerical methods to approximate the dichotomy spectrum).

For persistent synchronisation, we thus only require a dichotomy
spectrum in the directions transverse to the synchronisation manifold.
Instead we can impose the stricter condition of normal hyperbolicity
(see~\cite{Fenichel_71_1,Hirsch_77_1} and e.g.~\cite{Josic_00_1} in
the context of synchronisation of networks). That is, we also require
that any exponential contraction tangent to the synchronisation manifold
is weaker than in the transverse directions. In other words, the
spectra in the normal and tangential directions must be disjoint and
the normal spectrum must be strictly below the tangential one. In our
explicit setup, this so-called spectral gap condition translates to
\begin{equation*}
  \rho - \alpha \gamma < -r\,\rho \qquad\mbox{with $r \ge 1$.}
\end{equation*}
Under these assumptions we find a stronger form of persistence. Under
arbitrary $C^1$-small perturbations, solutions not only converge into
a neighbourhood of the synchronisation manifold, but an invariant
manifold\footnote{%
  Both smoothness and uniqueness of this manifold are subtle issues.
  In general the invariant manifold cannot be expected to be smoother
  than $C^r$. If the synchronisation manifold has a boundary (where it
  is only forward invariant), then non-uniqueness follows from local
  modifications that have to be made to apply the persistence theorem,
  see~\cite{Josic_00_1}. Note that both results hold, also when the
  synchronisation manifold $S$ is noncompact, see~\cite[Thm~3.1 and
  Chap.~4]{Eldering_13_1}.
}
\begin{equation*}
  \tilde{S} = \{ x_i = h_i(s), s \in U \subset \mathbb{R}^m, 1 \le i \le n \}
\end{equation*}
close to $S$ persists to which these solutions converge.
Moreover a stronger `shadowing' or `isochrony' property holds that any
solution curve $X(t)$ that converges to $\tilde{S}$, actually
converges at exponential rate $\tilde{\mu}$ to a unique solution curve
$X_{\tilde{S}}(t)$ on $\tilde{S}$ in the sense that there exists a $C$
such that for all $t \ge 0$
\begin{equation*}
  \norm{X(t) - X_{\tilde{S}}(t)} \le C e^{-\mu t}\,,
\end{equation*}
with $\mu$ close to $\alpha \gamma - \rho$.

\bigskip

\noindent \textbf{Acknowledgements.}
Tiago Pereira was supported by a Marie Curie IIF Fellowship (Project 303180), Jaap Eldering was supported by the ERC Advanced Grant 267382, and Martin Rasmussen and Jaap Eldering were supported by an EPSRC Career Acceleration Fellowship (2010--2015). We also thank CNPq and the Marie Curie IRSES staff exchange project DynEurBraz.

\newpage

\noindent
\textbf{References.}

\medskip


\end{document}